\documentclass[a4paper,12pt]{article}

\usepackage{amssymb,amsthm,amscd}
\usepackage{amsmath}
\usepackage{graphics}
\usepackage{latexsym,epic,eepic,epsfig}
\usepackage{color}

\usepackage{tocloft}
\usepackage{enumitem}

\theoremstyle{plain}
\newtheorem{theorem}{Theorem}
\newtheorem{defn}[theorem]{Definition}

\newtheorem{la}[theorem]{Lemma}
\newtheorem{cor}[theorem]{Corollary}
\newtheorem{rk}[theorem]{Remark}

\newcommand{\N}{\mathbb{N}}

\newcommand{\R}{\mathbb{R}}

\newcommand{\Z}{\mathbb{Z}}

\def\barr{\begin{array}}
\def\earr{\end{array}}
\def\beqarr*{\begin{eqnarray*}}
\def\eeqarr*{\end{eqnarray*}}

\def\mapright#1{\smash{\mathop{\longrightarrow}\limits^{#1}}}
\def\maprdown#1{\Big\downarrow
                   \rlap{$\vcenter{\hbox{$\scriptstyle#1$}}$}}
\def\mapldown#1{\llap{$\vcenter{\hbox{$\scriptstyle#1$}}$}\Big\downarrow}

\setitemize{noitemsep,leftmargin=15pt, topsep=0pt,parsep=0pt,partopsep=10pt}
\setenumerate{noitemsep,leftmargin=25pt, topsep=0pt,parsep=0pt,partopsep=0pt}
\begin{document}

\pagenumbering{arabic}


\title{\Large\bf  More on the Concept of Anti-integrability for H\'{e}non Maps}

\author{Zin Arai  \\
Department of Mathematical and Computing Science, School of Computing \\
Institute of Science Tokyo, Tokyo, Japan \\
Email: zin@comp.isct.ac.jp \\
 Yi-Chiuan Chen\\
   Institute of  Mathematics,  Academia Sinica, Taipei 106319, Taiwan\\
   Email: YCChen@math.sinica.edu.tw \\
 }

\date{8 May 2025}

\maketitle
\begin{center}
\end{center}


\begin{abstract}
  For the family of  H\'{e}non maps $(x,y)\mapsto (\sqrt{a}(1-x^2)-b y,x)$ of $\R^2$, the so-called anti-integrable (AI) limit concerns the limit $a\to\infty$ with fixed Jacobian $b$. At the AI limit, the dynamics reduces to a subshift of finite type.
There is a one-to-one correspondence between
sequences allowed by the subshift  and the AI orbits. The theory of anti-integrability says that each AI orbit can be continued to becoming a genuine orbit of  the H\'{e}non map for $a$ sufficiently large (and fixed Jacobian).

In this paper, we assume $b$ is a smooth function of $a$ and show that the theory can be extended to investigating the limit $\lim_{a\to\infty} b/\sqrt{a}=\hat{r}$ for any $\hat{r}>0$ provided that the one dimensional quadratic map $x\mapsto \displaystyle\frac{1}{\hat{r}}(1-x^2)$ is hyperbolic. 
\end{abstract}


\vspace{0.5cm}



\section{Introduction}

The concept of anti-integrability was first introduced by Aubry and Abramovici \cite{AA1990} in 1990 in order to prove the existence of chaotic trajectories in the standard map. 
At the anti-integrable (AI) limit, a dynamical system becomes non-deterministic, virtually a subshift of finite type. 
In the nice review \cite{Aubr1995}, Aubry gave more examples of dynamical systems that possess the AI limit. In particular, he showed that the singular limit $a\to \infty$ is an AI limit for the H\'{e}non map
\[ H_{a,b}:\mathbb{R}^2\to\mathbb{R}^2, \quad (x_i,y_i)\mapsto (x_{i+1}, y_{i+1})=(\sqrt{a}(1-x_i^2)-by_i,x_i), \quad i\in\Z,
\]
with fixed nonzero $b$, and that the trajectories at the limit can be characterised by the symbolic sequences of $\{-1,1\}^\mathbb{Z}$. In 1998, Sterling and Meiss \cite{SM1998} went further to show that
 every bounded orbit of $H_{a,b}$ is continued from a unique symbolic sequence belonging to $\{-1,1\}^\mathbb{Z}$ as long as $a$  complies with
\[
   a\ge \frac{(5+2\sqrt{5})(1+|b|)^2}{4} \qquad (\mbox{fixed}~b\not= 0).
\]
The parameter region above  was obtained first in \cite{DN1979} by Devaney and Nitecki and later in \cite{SM1998} using a different method (see also \cite{Chen2018}). 
Remark that the above parameter region is contained in the horseshoe locus and can be improved to be 
\[ a>2(1+|b|)^2, \quad b\not=0
\]
by taking the advantage of the Poincar\'{e} metric of complex analysis (see \cite{Ishi2008, MNTU2000, Mumm2008}). 

The H\'{e}non map $H_{a,b}$ has Jacobian $b$, and is invertible for non-vanishing  $b$. If $b=0$, it maps all of $\mathbb{R}^2$ to the curve given by $x=\sqrt{a}(1-y^2)$, thus collapses to the following quadratic map of one dimension:
\begin{equation}
 Q_{a}:\mathbb{R}\to\mathbb{R}, \quad x_i\mapsto x_{i+1}=\sqrt{a}(1-x_i^2), \quad i\ge 0.  \label{1Dquadratic}
\end{equation}
Therefore, it is clear that the concept of AI limit can also  be employed to study the  map $Q_{a}$, see \cite{Chen2007}. When $a>2$, it is well-known that the set of points which are bounded under the iteration of $Q_a$ is a Cantor set $\Lambda^\prime_a$ of measure zero, and the restriction of $Q_a$ to $\Lambda^\prime_a$ is an embedded one-sided full shift with two symbols. Hence, from the AI limit point of view, the embedded shift persists from the limit $a\to\infty$ to $a=2^+$. Similarly, the singular limit $a\to \infty$ is also an AI limit for a mapping of the form
\[  x_i\mapsto x_{i-1}=\sqrt{a}(1-x_i^2), \quad x_i\in\mathbb{R},~ i\le 0.
\]

Let $\epsilon=1/\sqrt{a}$, $r=b/\sqrt{a}$. Apparently, for nonzero $a$, $b$, $\epsilon$ and $r$, a sequence $\left(x_i,x_{i-1}\right)_{i\in\mathbb{Z}}$ is an orbit of $H_{1/\epsilon^2,r/\epsilon}$ if and only if the sequence $(x_i)_{i\in\mathbb{Z}}$ is a solution of the following second order recurrence relation
\begin{equation}
  \epsilon x_{i+1}-(1-x_i^2)+r x_{i-1}=0, \quad i\in\mathbb{Z}. \label{HenonRecurrence}
\end{equation}
The AI limit established in \cite{Aubr1995, SM1998} for the H\'{e}non map can be viewed as the limiting situation $(\epsilon,r)\to (0,0)$ with $r=\epsilon b$ for some nonzero constant $b$, whereas the AI limit revealed in \cite{Chen2007} is the  situation $(\epsilon,r)\to (0,0)$ with fixed  $r=0$. (Note that, when $\epsilon= r=0$, equation \eqref{HenonRecurrence} becomes an algebraic one, and can be solved easily to obtain $x_i=-1$ or $1$ for each $i$; whether $x_i$ is equal to $-1$ or $1$ does not determine the value of $x_{i+1}$, therefore the system  is non-deterministic; for any $x_i$ in a given $(x_i)_{i\in\Z}\in\{-1,1\}^\Z$, the value of $x_{i+1}$ can be obtained by the usual shift operation $(x_i)_{i\in\Z}\mapsto (\tilde{x}_i)_{i\in\Z}$ with $\tilde{x}_i=x_{i+1}$, therefore the system is a Bernoulli shift with two symbols. Hence, $(\epsilon, r)=(0,0)$ is an AI limit for \eqref{HenonRecurrence} in the sense of \cite{AA1990}.)

\begin{rk} \rm
Another popular form of H\'{e}non map is $\mathcal{H}_{a,b}:(x,y)\mapsto (a-x^2+by, x)$.    Let $\mathcal{R}=\mathcal{R}(a,b):=\frac{1}{2}\left(1+|b|+\sqrt{(1+|b|)^2+4a}\right)$. Devany and Nitecki \cite{DN1979} showed that the bounded orbits of  $\mathcal{H}_{a,b}$ are confined in the region $\left\{   (x,y)\in\R^2|~\max \{|x|, |y|\}\le \mathcal{R}\right\}$.
The homeomorphism $(x,y)\mapsto (\sqrt{a}x, \sqrt{a}y)$ conjugates $\mathcal{H}_{a,b}$ to  $H_{a,-b}$  provided that $a$  is positive. 
 Therefore, the bounded orbits of $H_{1/\epsilon^2, r/\epsilon}$ are confined in the region  $\left\{   (x,y)\in\R^2|~\max \{|x|, |y|\}\le \frac{1}{2}\left(\epsilon+|r|+\sqrt{(\epsilon+|r|)^2+4}\right)\right\}$.
\end{rk}

As has been pointed out in the seminal paper \cite{AA1990} that  a dynamical system may have more than one AI limit. An example is   the following two-harmonic area-preserving twist map of $\mathbb{R}^2$ with two parameters $k_1$ and $k_2$ defined by
\begin{subequations} \label{2harmonicMap}
\begin{eqnarray}
  x_{i+1} &=& x_i+y_{i+1} \quad (\mbox{mod 1}), \\
  y_{i+1} &=& y_i-\frac{k_1}{2\pi}\sin 2\pi x_i-\frac{k_2}{4\pi}\sin 4\pi x_i, \qquad i\in\mathbb{Z}.
\end{eqnarray}
\end{subequations}
If we represent the space of parameters $(k_1,k_2)$ as the extended complex plane, then the limit $k_1+ik_2\to \infty$ along a path $k_2/k_1=k\in\mathbb{R}$ may resulting in a different AI limit for a different choice of $k$, see Baesens and MacKay \cite{BM1993}.  

With these in mind, let us return back to the H\'{e}non family $H_{a,b}$, which also has two parameters, $a$ and $b$. We have seen that the H\'{e}non family possesses the AI limits $(a,b)\to\infty$, viewed in the extended complex $(a+ib)$-plane, along the path $b=nonzero~ constant$, for which the limiting dynamics $\lim_{a\to\infty}H_{a,b}$ reduces to the two-sided full shift with two symbols, and along the path $b=0$, for which $\lim_{a\to\infty}H_{a,0}$  reduces  to the one-sided full shift with two symbols (more precisely, to the two-sided shift on the inverse limit space of one-sided full shift). It is fairly natural to ask what the limiting dynamics of the H\'{e}non map would be if the limit $(a,b)\to\infty$ is taken along an arbitrary path, for example along the path $b=\log a$ or $b=\sqrt[3]{a}$ or the path $b/\sqrt{a}=constant$.
 We believe that a  study of this question can shed more light on the  concept of   anti-integrability and on its methodology that has been developing to prove the existence of chaotic orbits in a variety of systems. 

\begin{rk} \rm
~{}
\begin{itemize}
\item
In fact, we get $(\epsilon, r)\to (0,0)$ in the recurrence relation (\ref{HenonRecurrence}) when $(a,b)\to\infty$ along a path $b/\log a=constant$ or $b/\sqrt[3]{a}=constant$, resulting in an AI limit.
\item
The roles played by $\epsilon$ and $r$ in  (\ref{HenonRecurrence}) are symmetric in the sense that exchanging $\epsilon$ and $r$ is equivalent to reversing the time. (This also shows that the inverse $H_{a,b}^{-1}(x,y)$ is topologically conjugate to $H_{a/b^2,1/b}(x,y)$ by exchanging $x$ and $y$.)  Hence,  rather than the path $b/a=constant$, the path $b/\sqrt{a}=constant$ is a more  natural one to look at.
\end{itemize}
\end{rk}
\begin{rk} \rm
Apart from the standard, the H\'{e}non and the quadratic maps, the concept  of AI limit 
has been applied to various  dynamical systems, such as   high-dimensional symplectic maps \cite{Chen2005, MM1992} and H\'{e}non-like maps \cite{Qin2001}, as well as the Smale horseshoe \cite{Chen2006}. This concept and methodology have shown powerful for proving existence of chaotic invariant sets not only in discrete-time systems  but also in continuous-time ones, for instance,   Hamiltonian and Lagrangian systems \cite{BCM2013, BM1997, Chen2003}, $N$-body problems \cite{BM2000}, and  billiard systems \cite{Chen2004, Chen2010}.  See \cite{BT2015} for a nice review. 
\end{rk}

The structure of the rest of this paper is as follows. In the next section, we present our main result in Theorem \ref{thm:AI2}, which 
treats parameter $b$ as a function of parameter $a$, and can be regarded as a generalisation of  the theory of anti-integrability for H\'enon maps. 
  The classical result of the theory treated in \cite{Aubr1995, SM1998}  is  when $b$ equals a constant. We shall see in Theorem \ref{thm:AI2}  that all statements of the theory of anti-integrability  are also valid when the limit $\lim_{a\to\infty}a/b^2$ exists and belongs to a hyperbolic parameter of the quadratic map \eqref{1Dquadratic}.  In other words, the theory of anti-integrability can be extended to investigating the limit $\lim_{a\to \infty}b/\sqrt{a}=\hat{r}$ for any $\hat{r}>0$ provided that $x\mapsto (1-x^2)/\hat{r}$ is uniformly hyperbolic. Our approach to Theorem \ref{thm:AI2} is via the implicit function theorem. 
A peculiar feature of our approach is that it automatically gives rise to the so-called {\it quasi-hyperbolicity} in the sense of \cite{CFS1977, SS1974}. 
 We show  in Section \ref{sec:hyp} that the collection of orbits continued from the AI orbits of the H\'enon map forms a quasi-hyperbolic invariant set. 
The quasi-hyperbolicity  naturally links the persistence of uniformly hyperbolic invariant sets for parameters near $\epsilon=0$ or $r=0$ to the uniformly hyperbolic plateaus studied by the first author \cite{Arai2007}, since
quasi-hyperbolicity plus the chain recurrence of the set establish the uniform hyperbolicity \cite{CFS1977, SS1974}.
 Plateaus for varies interesting parameter spaces for the H\'enon maps are demonstrated also in Section 
\ref{sec:hyp}. In particular, any path for $(a,b)$ to the infinity may be visualised through figures in that section. Section \ref{sec:proof} is devoted to the
 proof of the main theorem (Theorem \ref{thm:AI2}) of this paper.

\section{Main results}

As mentioned in the preceding section, equation \eqref{HenonRecurrence} reduces to a one-dimensional quadratic map $x_i\mapsto x_{i-1}=(1-x_i^2)/r$ when $\epsilon$ tends to zero for a fixed nonzero $r$. Hence, the quadratic map $Q_a$ plays a crucial role in  the study of the recurrence relation \eqref{HenonRecurrence} with small $\epsilon$. 

It is well known that \cite{Sing1978} the quadratic map $Q_a$ has at most one attracting periodic orbit, whose immediate basin of attraction contains the critical point $x=0$. If $Q_a$ has an attracting periodic orbit or if the critical point goes to infinity under the iteration of $Q_a$, we say that the map $Q_a$ is {\it  hyperbolic}. Graczyk and \'Swiatek \cite{GS1997} showed that the set of parameter $a$ for which $Q_a$ is hyperbolic  is open and dense in the interval $[0,2]$. If $Q_a$ has an attracting periodic orbit, the complement of the basin of attraction of the  orbit is compact and invariant in both forward and backward iterations. It has also been known  that this complement is a uniformly hyperbolic set $\Lambda^{\prime\prime}_a$ of measure zero, and the restriction of $Q_a$ to $\Lambda^{\prime\prime}_a$ is topologically conjugate to a one-sided topological Markov chain, see e.g. \cite{KH1995, Mane1985, dMvS1993, Misi1981}. 

We use $\sigma$ to denote the usual left-shift map on the product space $\R^\Z$ or $\R^{\N_0}$, where $\mathbb{N}_0=\{0,1, 2, \ldots\}$.

\begin{defn} \rm
~{}
\begin{enumerate}
\item When $Q_{a}$ is hyperbolic, define 
\[ \Lambda_{a}:=\begin{cases}
    \Lambda^\prime_{a} & \mbox{for}~a\in(2, \infty) \\
\Lambda^{\prime\prime}_{a} & \mbox{for}~a \in (0, 2).
                         \end{cases}
\]
We use 
$\Lambda_{\infty}:= \{-1, 1\} 
$
and $Q_{\infty}:=\sigma|{\{-1,1\}^{\mathbb{N}_0}}$.
\item Define $\mathsf{Hyp}:=\{a|~ Q_a~\mbox{is hyperbolic}\}\cup\{\infty\}$. 
\end{enumerate}
\end{defn}

Now, we specify how parameters $(a,b)\to\infty$, viewed in the extended complex $(a+ib)$-plane,  by 
letting  
\[ r=r(\epsilon)=\epsilon b=b/\sqrt{a}
\]
 with $r$ depends continuously differentiable on $\epsilon$. (Note that we parameterise $b$ by $a$, but $b$ may be a constant.) Since the H\'{e}non map is invertible only for non-vanishing Jacobian, we assume $r(\epsilon)\not=0$ when $\epsilon\not=0$. Let
\[  \lim_{\epsilon\to 0}r(\epsilon)=:\hat{r}.
\]

Let $l_\infty(\mathbb{Z}, \mathbb{R}):=\{{\bf x}=(x_i)_{i\in\mathbb{Z}}|~ x_i\in\mathbb{R},~\mbox{bounded}~ \forall i \}$ be the Banach space of bounded bi-infinite sequences endowed with the supremum norm $\|{\bf x}\|=\sup_{i\in\Z}|x_i|$.
Then, a bounded sequence ${\bf x}$ that solves the recurrence relation \eqref{HenonRecurrence} precisely is a zero of the map $F(\cdot; \epsilon):l_\infty(\mathbb{Z},\mathbb{R})\to l_\infty(\mathbb{Z},\mathbb{R})$ defined by
\[ F({\bf x}; \epsilon)=\left( F_i({\bf x}; \epsilon)\right)_{i\in\Z},
\]
with 
\[
 F_i({\bf x}; \epsilon)=\epsilon x_{i+1}-(1-x_i^2)+r x_{i-1}.  
\]
Define subspaces of $l_\infty(\Z, \R)$,
\[
 \Sigma_{\epsilon} := \left\{{\bf x}\in l_\infty (\Z, \R)|~ F({\bf x}; \epsilon)=0\right\},
\]
and 
\[
 \Sigma_{0, \hat{r}} := \left\{{\bf x}\in l_\infty(\Z, \R)|~ F({\bf x}; 0)=0~\mbox{and}~x_i\in \Lambda_{1/\hat{r}^2} ~\forall i\in\Z\right\}
\]
when $1/\hat{r}^2\in\mathsf{Hyp}$.
Notice that 
\[   \Sigma_{0,0}=\{-1, 1\}^\Z,
\] 
and that 
\[ \Sigma_{0, \hat{r}}~\mbox{is homeomorphic to the inverse limit space}~ \underleftarrow{\lim}(\Lambda_{1/\hat{r}^2}, Q_{1/\hat{r}^2})
\]
 with the product topologies when $1/\hat{r}^2\in\mathsf{Hyp}$. Recall that
   the  {\it inverse limit space}  for a  continuous map $f$ of $\R^k$, $k\ge 1$, is defined by 
\[   \underleftarrow{\lim} (\R^k,f):=\{(\ldots,z_{-1},z_0)\in (\R^k)^{\mathbb{Z}^-}|~ f(z_{i-1})=z_i  ~ \forall ~ i\le 0\}.
\]
It induces a  map, which is often referred  to as the shift homeomorphism on 
$\underleftarrow{\lim} (\R^k,f)$,
\[
   \widehat{f}(\ldots,z_{-2},z_{-1},z_0) = (\ldots,z_{-1},z_0,f(z_0)).
\]
The map $f$ is called the {\it bonding map} of the inverse limit space.
Let $S$ be an $f$-invariant subset of $\R^k$, we use $
   \underleftarrow{\lim} (S,f):=\{(\ldots,z_{-1},z_0)
 \in\underleftarrow{\lim} (\R^k,f)|~ z_i\in S ~ \forall ~ i\le 0\}.$
 
\begin{defn} \rm
 A zero  ${\bf x}^*$ of $F(\cdot; \epsilon)$ is said to be {\it non-degenerate} if the linear map $DF({\bf x}^*; \epsilon): l_\infty(\Z,\R)\to l_\infty(\Z,\R)$, which is the derivative of $F(\cdot; \epsilon)$ at ${\bf x}^*$, is invertible. 
\end{defn}

Define the following set of $\R^2$, 
\[ 
      \mathcal{A}_{0}:=\bigcup_{{\bf x}\in \Sigma_{0,\hat{r}}}(x_0, x_{-1}).  
\]
The theorem below  is the main aim of this paper.

\begin{theorem} \label{thm:AI2}
For the family of H\'{e}non maps $H_{1/\epsilon^2, r/\epsilon}$   with $\lim_{\epsilon\to 0}r(\epsilon)=\hat{r}$ and $1/\hat{r}^2\in\mathsf{Hyp}$, we have
\begin{enumerate}
\item 
${\bf x}^\dag$ is a non-degenerate zero of $F(\cdot, 0)$ for  every ${\bf x}^\dag\in \Sigma_{0,\hat{r}}$. Hence, for every ${\bf x}^\dag\in\Sigma_{0,\hat{r}}$, there exist positive constants $\epsilon_0, \delta_0$, and a unique ${\bf x}^*={\bf x}^*(\epsilon;{\bf x}^\dag)\in\Sigma_{\epsilon}$ such that $|x_i^*-x_i^\dag|<\delta_0$ for all  $i\in\Z$ provided $0\le \epsilon<\epsilon_0$.
\item
The constants $\epsilon_0$, $\delta_0$ are independent of ${\bf x}^\dag$.
 Let $\mathcal{A}_{\epsilon}:=\bigcup_{{\bf x}^\dag\in \Sigma_{0,\hat{r}}}(x_0^*, x_{-1}^*)$. The restriction of  $H_{1/\epsilon^2, r/\epsilon}$  to the compact invariant set $\mathcal{A}_{\epsilon}$ is topologically conjugate to $\sigma$ on $\Sigma_{0,\hat{r}}$  provided $0<\epsilon<\epsilon_0$.
\item $x_i^*\to x_i^\dag$ for all $i$ as $\epsilon\to 0$,  and $\mathcal{A}_{\epsilon}\to \mathcal{A}_{0}$ in the Hausdorff topology.    
\end{enumerate}
\end{theorem}

The proof of  Theorem \ref{thm:AI2}  is presented in Section \ref{sec:proof}.
The theorem  tells that the dynamics of the H\'{e}non map $H_{1/\epsilon^2, r/\epsilon}$ on $\mathcal{A}_{\epsilon}$ collapses to that of $\sigma$ on $\Sigma_{0, \hat{r}}$ as $\epsilon$ tends to zero. 
Moreover, in terms of parameters $a$ and $b$,
  Theorem \ref{thm:AI2} implies  that the H\'{e}non map $H_{a,b}$ restricted to the set $\mathcal{A}_{1/\sqrt{a}}$ is topologically conjugate to $\sigma$ on $\Sigma_{0,\hat{r}}$ if $1/\epsilon_0^2<a<\infty$. 

The theory of anti-integrability  for H\'enon maps is a special case of Theorem \ref{thm:AI2}. For a fixed $b$, setting $r(\epsilon)=\epsilon b$ (thus $\lim_{\epsilon\to 0}r(\epsilon)=\hat{r}=0$), we arrive at
\begin{cor}[Anti-integrability for H\'enon maps] \label{thm:AI}
 The family of H\'{e}non maps $H_{1/\epsilon^2, b}$ 
has an AI  limit $\epsilon\to 0$   for  any  fixed $0<b<\infty$:
\begin{enumerate}
\item 
${\bf x}^\dag$ is a non-degenerate zero of $F(\cdot, 0)$ for  every ${\bf x}^\dag\in \Sigma_{0,0}=\{-1,1\}^\Z$. Hence, for every ${\bf x}^\dag\in\{ -1, 1\}^{\Z}$, there exist positive constants $\epsilon_0, \delta_0$, and a unique ${\bf x}^*={\bf x}^*(\epsilon;{\bf x}^\dag)\in\Sigma_{\epsilon}$ such that $|x_i^*-x_i^\dag|<\delta_0$ for all  $i\in\Z$ provided $0\le \epsilon<\epsilon_0$.
\item
  The restriction of  $H_{a, b}$  to  $\mathcal{A}_{1/\sqrt{a}}$ is topologically conjugate to   $\sigma$ on $\{-1, 1\}^\Z$  provided that $1/\epsilon_0^2<a<\infty$. 
\end{enumerate}
\end{cor}

Now, let us consider the case $b/\sqrt{a}$ equal to a fixed constant $\hat{r}$, as mentioned in the Introduction section.
Since $\Sigma_{0,\hat{r}}$ is homeomorphic to $\underleftarrow{\lim}(\Lambda_{1/\hat{r}^2}, Q_{1/\hat{r}^2})$, by letting $a/b^2=1/\hat{r}^2$ we have  

 \begin{cor}[$b/\sqrt{a}=$ a fixed nonzero constant]
Suppose that $a/b^2$ is a constant  and belongs to $\mathsf{Hyp}$. The restriction of  $H_{a, b}$  to  $\mathcal{A}_{1/\sqrt{a}}$ is topologically conjugate to    $\widehat{Q}_{a/b^2}$ on $\underleftarrow{\lim}(\Lambda_{a/b^2}, Q_{a/b^2})$   provided that $1/\epsilon_0^2<a<\infty$.  
 \end{cor}

Concerning the topological entropy $h_{\rm top}$, it is well-known \cite{Bowe1970} that the topological entropy of the shift homeomorphism of the inverse limit space is equal to the one of the bonding map. In other words, $h_{\rm top}(\widehat{Q}_{a/b^2}|\underleftarrow{\lim}
(\Lambda_{a/b^2}, Q_{a/b^2}))=h_{\rm top}(Q_{a/b^2}|\Lambda_{a/b^2})$. As a consequence, 
\begin{cor} \label{cor:AI_Henon}
For $1/\epsilon_0^2<a<\infty$, 
\begin{enumerate}
\item
$h_{\rm top}(H_{a,b}|\mathcal{A}_{1/\sqrt{a}}) = h_{\rm top}(\sigma| \{-1,1\}^\Z) =h_{\rm top}(\sigma| \{-1,1\}^{\N_0})=\log 2$ if $b$ is a fixed nonzero constant. 
\item
$h_{\rm top}(H_{a,b}|\mathcal{A}_{1/\sqrt{a}}) = h_{\rm top}(Q_{a/b^2}| \Lambda_{a/b^2}) =h_{\rm top}(Q_{a/b^2})$ 
if   $b/\sqrt{a}$ equals  a  nonzero constant and  $Q_{a/b^2}$ is hyperbolic. 
\end{enumerate}
\end{cor}

\section{Quasi-hyperbolicity and hyperbolic plateaus} \label{sec:hyp}

\begin{defn} \rm
An invariant set  $S$  of  a $C^1$ diffeomorphism $f$ of $\R^k$, $k\ge 2$, is called {\it quasi-hyperbolic}  if the recurrence relation
 \begin{equation}
   \zeta_{i+1}-Df(z_i)~ \zeta_i=0,\quad i\in\mathbb{Z}, \label{hypequivhomo}
 \end{equation}
 has no non-trivial bounded solutions  $(\zeta_i)_{i\in\Z}$ for any orbit $(z_i)_{i\in\Z}$ with $z_0\in S$. 
\end{defn}

\begin{cor} \label{cor:quasi}
Providing that $\lim_{\epsilon\to 0} r(\epsilon)=\hat{r}$, $1/{\hat r}^2\in \mathsf{Hyp}$, and $0<\epsilon<\epsilon_0$, the set $\mathcal{A}_{\epsilon}$  is quasi-hyperbolic for $H_{1/\epsilon^2, r/\epsilon}$.
\end{cor}
\proof
Since $\mathcal{A}_{\epsilon}=\bigcup_{{\bf x}^\dag\in \Sigma_{0,\hat{r}}}(x_0^*, x_{-1}^*)$ with ${\bf x}^*={\bf x}^*(\epsilon; {\bf x}^\dag)\in \Sigma_{\epsilon}$, the operator $DF({\bf x}^*;  \epsilon)$ is invertible provided that $0\le \epsilon<\epsilon_0$. 
When $\epsilon\not=0$, the invertibility implies that the recurrence relation
 \[
  \epsilon\xi_{i+1}+2 x_i \xi_i + r \xi_{i-1} = \epsilon\eta_i-r\tilde{\eta}_{i-1}
  \] has a unique solution $(\xi_i)_{i\in\mathbb{Z}}\in l_\infty(\Z, \R)$ for any given $(\eta_i)_{i\in\Z}$, $(\tilde{\eta}_i)_{i\in\Z}\in l_\infty(\Z,\R)$.
Let  $\tilde{\xi}_{i}=\xi_{i-1}+\tilde{\eta}_{i-1}$. Then, the following recurrence relations 
 \begin{eqnarray*}
  \xi_{i+1}+\frac{2 x_i}{\epsilon} \xi_i +\frac{r}{\epsilon}\tilde{\xi_i} &=& \eta_i, \\
  \tilde{\xi}_{i+1}-\xi_i &=&\tilde{\eta}_i,
 \end{eqnarray*}
or equivalently, the one
\[ 
 \left ( \begin{matrix}
                                           \xi_{i+1}  \\
                                         \tilde{\xi}_{i+1}  \end{matrix}
                              \right ) -
\left (\begin{matrix} 
        -\displaystyle\frac{2 x_i}{\epsilon}  &   -\displaystyle\frac{r}{\epsilon}\\
1 &0 
\end{matrix} \right ) \left ( \begin{matrix}
                                           \xi_i  \\
                                         \tilde{\xi}_i  \end{matrix}
                              \right ) = \left ( \begin{matrix}
                                           \eta_i  \\
                                         \tilde{\eta}_i  \end{matrix}
                              \right )
\]
 has a unique  solution  $(\xi_i,\tilde{\xi}_{i})_{i\in\mathbb{Z}}\in  l_\infty(\Z,\R^2)$ for any given $(\eta_i, \tilde{\eta}_i)_{i\in\Z}\in l_\infty(\Z,\R^2)$. This further implies that the only solution  $(\xi_i,\tilde{\xi}_{i})_{i\in\mathbb{Z}}$ lying in the space $l_\infty(\Z,\R^2)$ for 
 \[ 
\left ( \begin{matrix}
                                           \xi_{i+1}  \\
                                         \tilde{\xi}_{i+1}  \end{matrix}
                              \right ) -
\left (\begin{matrix} 
        -\displaystyle\frac{2 x_i}{\epsilon}  &  -\displaystyle\frac{r}{\epsilon}\\
  1&0 
\end{matrix} \right ) \left ( \begin{matrix}
                                           \xi_i  \\
                                         \tilde{\xi}_i  \end{matrix}
                              \right ) = 0
\]
is zero. Since the two-by-two matrix in the above equation is just $DH_{1/\epsilon^2, r/\epsilon}(x_i, y_i)$ with $y_i=x_{i-1}$ and $(x_0, y_0)\in \mathcal{A}_{\epsilon}$, the set $\mathcal{A}_{\epsilon}$ is quasi-hyperbolic. 
 \qed

Denote by $\mathcal{CR}(f)$ the chain-recurrent set of a map $f$. For a compact $f$-invariant set $S$, if it is chain-recurrent, namely $\mathcal{CR}(f|S)=S$, 
Churchill {\it et al}. \cite{CFS1977}, Sacker and Sell \cite{SS1974} showed the following result:
\begin{theorem} \label{thm:quasi}
A  compact chain-recurrent quasi-hyperbolic invariant set $S$ of  a $C^1$ diffeomorphism $f$ of $\R^k$, $k\ge 2$, is uniformly hyperbolic.
\end{theorem}

Corollary \ref{cor:quasi} and Theorem \ref{thm:quasi} yield the following:
\begin{cor} \label{cor:unihyp}
 Let $J_{\epsilon}$ be a compact chain-recurrent invariant subset of  $\mathcal{A}_{\epsilon}$, namely $\mathcal{CR}(H_{1/\epsilon^2, r/\epsilon}|J_{\epsilon})=J_{\epsilon}$, then $J_{\epsilon}$ is uniformly hyperbolic.
\end{cor}

Corollary \ref{cor:unihyp} is closely related to the work of \cite{Arai2007}. There, based on a rigorous numerical  method, the first author of the current paper obtained parameter regions, called the {\it uniformly hyperbolic plateaus} (see Figure 1 of \cite{Arai2007}), on the $(a,b)$-plane for which the H\'{e}non map  $\mathcal{H}_{a,b}$ is uniformly hyperbolic on its chain-recurrent set.

In the remainder of this section, we embed these plateaus   into varies parameter spaces. On one hand, seeing how the plateaus are presented  is interesting in its own right. Especially, because the parameter region $\{(a,b)|~a>0\}$ is unbounded, it is desirable to make a transformation of parameters so that the ``whole" picture of  the hyperbolic plateaus can be seen in a bounded domain.
On the other hand, these plateaus provide vivid pictures of how the limiting dynamics of the H\'{e}non map would be as parameters approach some limit along a certain path.

\begin{figure}[htbp]
	\begin{center}
\includegraphics[width=0.40\textwidth]{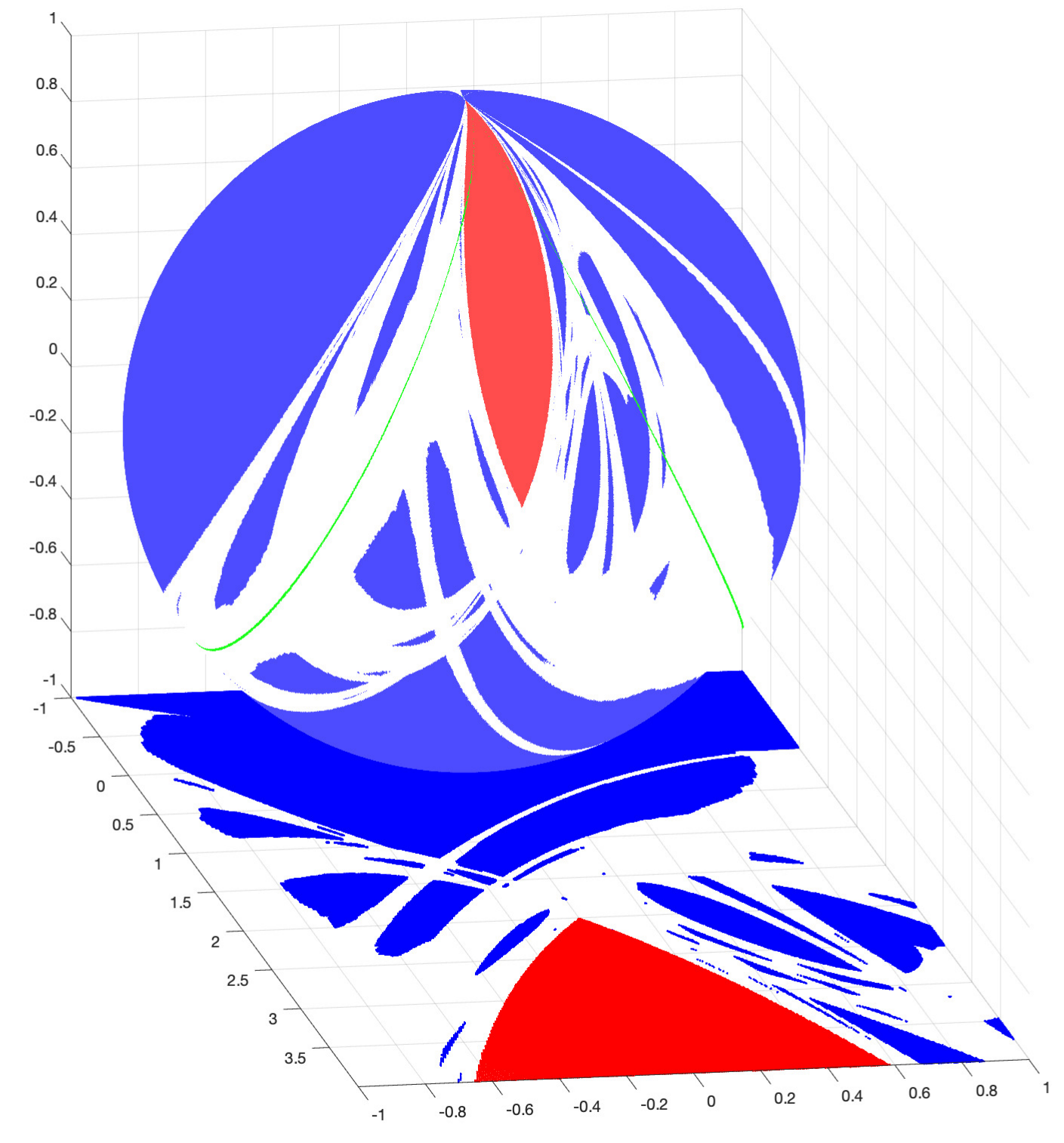}
\includegraphics[width=0.58\textwidth]{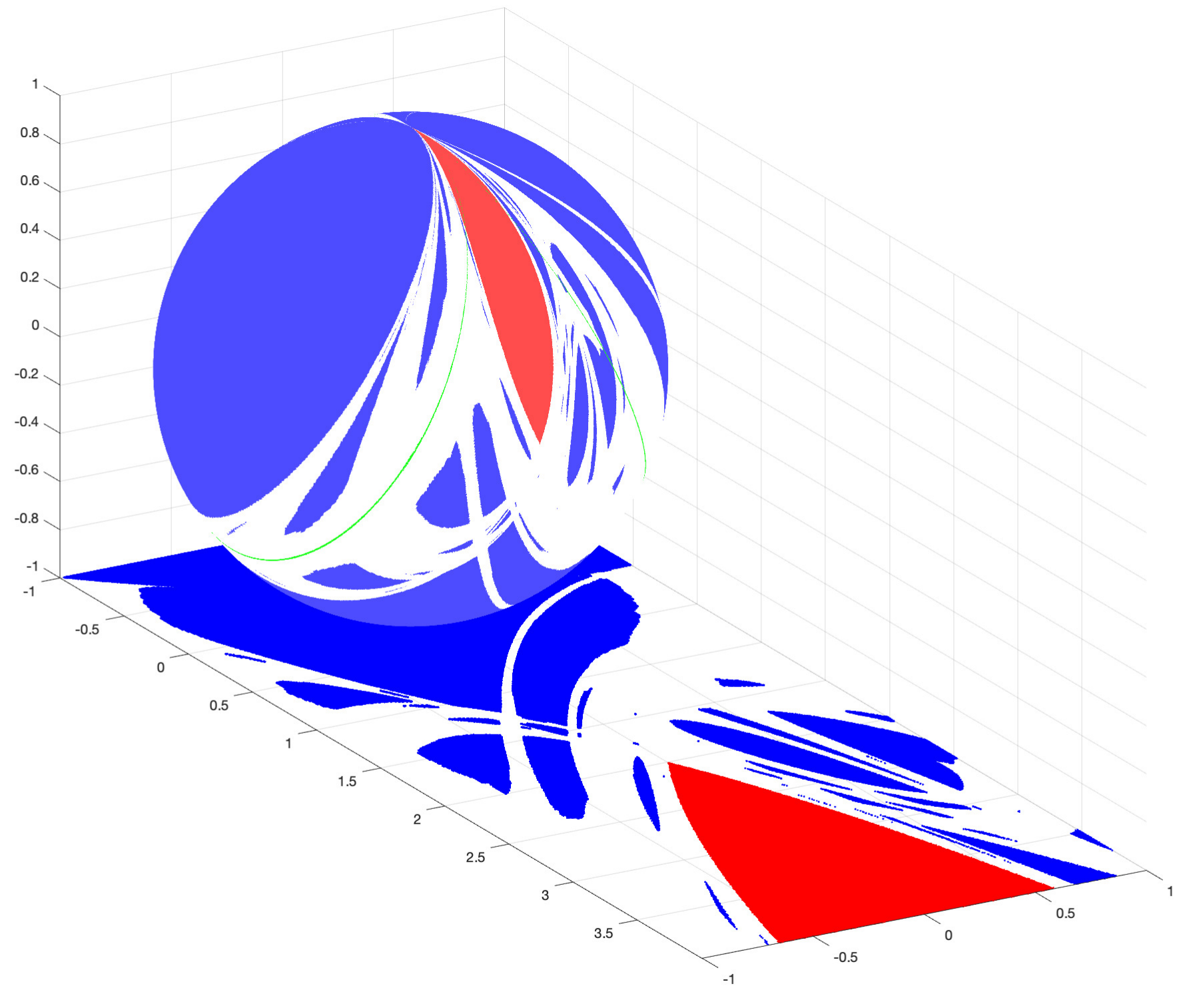}	
\end{center}
	\caption{Hyperbolic plateaus on the sphere and its  stereographic projection viewed from two different angles.}
	\label{fig:plateaus_sphere}
\end{figure}

\paragraph{Plateaus on sphere.}
Using the stereographic projection, we can convert the plateaus from the $(a,b)$-plane onto the unit sphere $\{(X,Y,Z)\in\R^3|~X^2+Y^2+Z^2=1\}$ by $(X, Y, Z)=(4a/(a^2+b^2+4), 4b/(a^2+b^2+4), (a^2+b^2-4)/(a^2+b^2+4))$. See Figure \ref{fig:plateaus_sphere}. Note that the projection maps the ``south pole" $(X,Y,Z)=(0,0,-1)$ to the origin $(a,b)=(0,0)$, while the ``north pole" to the ``infinity" $a^2+b^2=\infty$. 
The coloured  part  in the figure  is the parameter region 
where  the uniform hyperbolicity of the chain recurrent set of the H\'{e}non map is proved numerically. The region depicted in red is contained in the horseshoe locus. The green curves on the sphere under the projection are the two lines  $b=\pm 1$ on the $(a,b)$-plane.
As parameters $(a,b)\to\infty$ along a given path,  the limiting dynamics of  the H\'{e}non map $H_{a,b}$  would depend on how the corresponding curve approaches the north pole on the sphere. 

\begin{figure}[htbp]
	\begin{center}
		\includegraphics[width=0.490\textwidth]{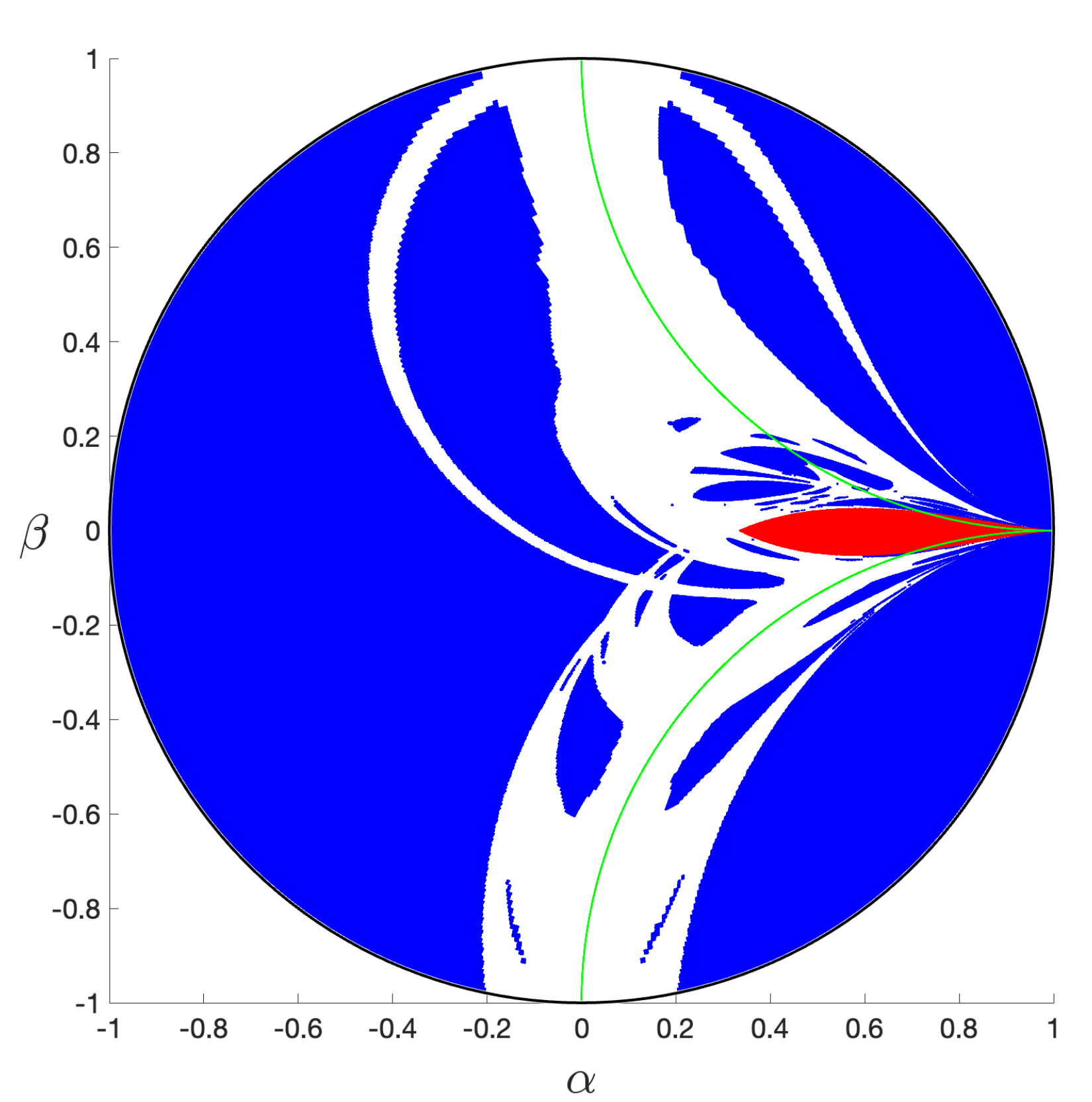}
\includegraphics[width=0.490\textwidth]{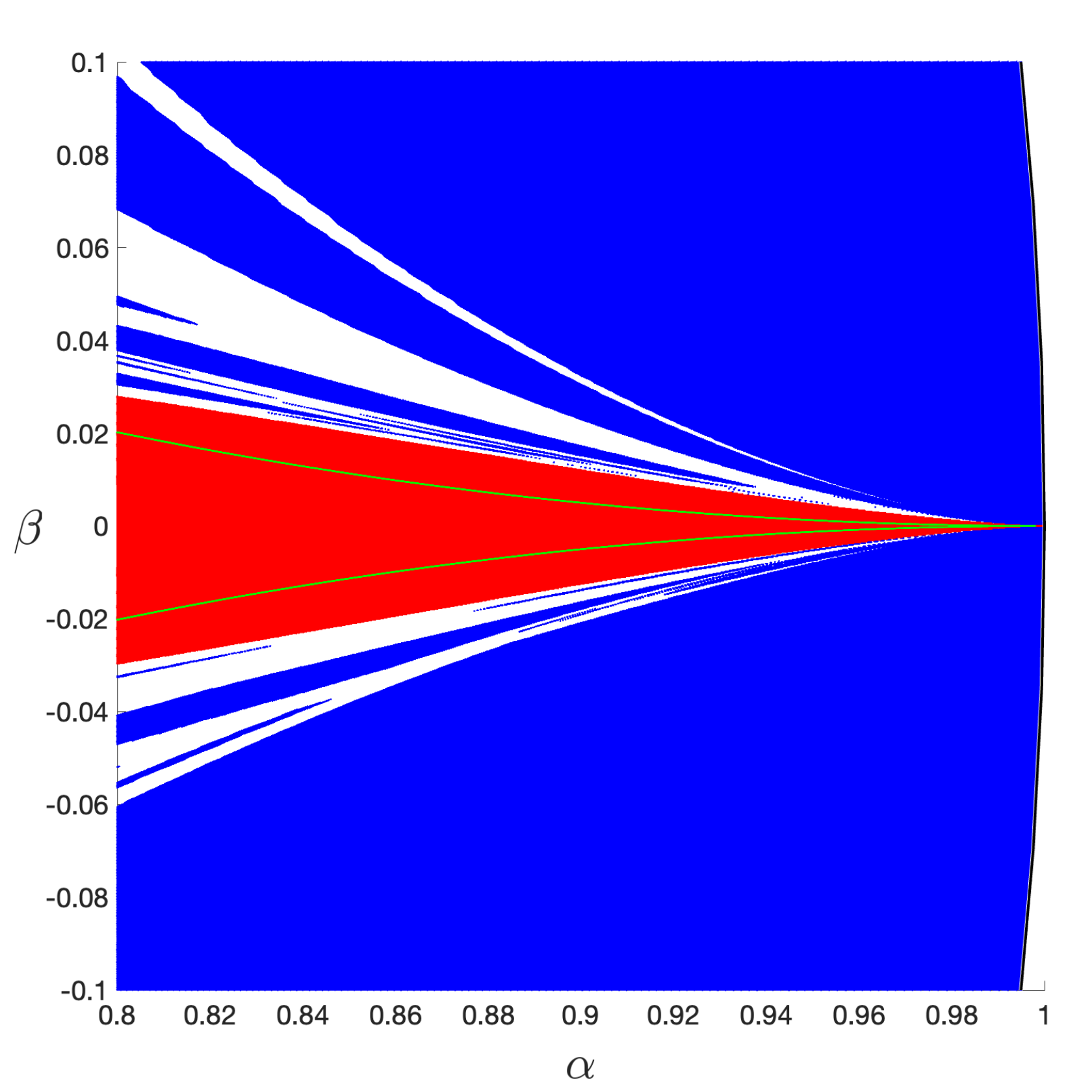}
	\end{center}
	\caption{(Left) Hyperbolic plateaus  under  the M\"{o}bius transformation $\mathcal{M}$.  (Right) Magnification for $0.8\le \alpha<1$ and $-0.1\le \beta \le 0.1$.}
	\label{fig:plateaus_disc}
\end{figure}

\paragraph{Plateaus on  disc.}
Since the parameter $a$ for $H_{a,b}$ must be non-negative, the
 following M\"{o}bius transformation  $\mathcal{M}: a+ib\mapsto (a+ib-1)/(a+ib+1)=:\alpha+i\beta$, which maps the right half complex plane $\{a+ib|~a>0\}$ to the open unit disc $\{\alpha+i\beta|~0\le \alpha^2+\beta^2<1\}$, transforms  the plateaus in the $(a,b)$-plane to plateaus in the $(\alpha, \beta)$-plane. The result is shown in Figure \ref{fig:plateaus_disc}.
Again, the coloured regions in the figure are the hyperbolic plateaus.  The image of the two horizontal lines $\{a+ib|~a>0,~ b=\pm 1\}$ under  $\mathcal{M}$ is depicted in green colour. The horseshoe locus contains the red region.
 Notice that the M\"{o}bius transformation preserves angles, and maps  lines to lines or circles, and circles to  lines or circles.

\begin{figure}[htbp]
	\begin{center}
		\includegraphics[width=0.33\textwidth, angle=0]{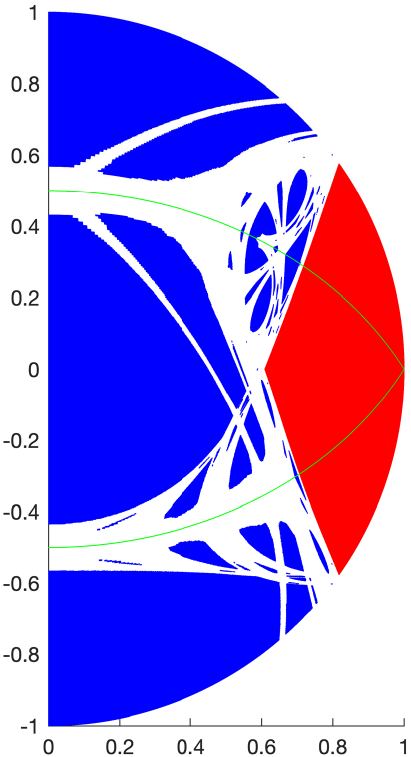}
	\end{center}
	\caption{Hyperbolic plateaus in the $(\rho,\theta)$-coordinates.}
	\label{semi_disc}
\end{figure}

\paragraph{Plateaus on semi-disc.}
 Let $\rho=(2/\pi)\tan^{-1}\sqrt{a+b^2}$ and $\theta=\tan^{-1}(b/\sqrt{a})$. Figure \ref{semi_disc} shows the plateaus on a semi-disc of radius one in the $(\rho, \theta)$-coordinates. The semi-circle $\{(\rho, \theta)|~ \rho=1,~ -\pi/2<\theta<\pi/2\}$ is where   the ``infinity" $(a,b)=(\infty, b)$ is.  The coloured regions are the hyperbolic plateaus. The red region belongs to the horseshoe locus. The green curves are depicted by $\{(2/\pi)\tan^{-1}\sqrt{a+1}, \tan^{-1}(\pm 1/\sqrt{a}))|~ 0<a<\infty\}.$
An advantage of displaying the plateaus in this polar-coordinates is that if parameters $a$ and $b$ tend to infinity with a distinct limit $\lim_{a\to\infty}b/\sqrt{a}=\hat{r}$,
 then the corresponding path on the semi-disc will approach  a distinct point $(1, \tan^{-1}\hat{r})$. Figure \ref{semi_disc} may be compared with Figure 2 of \cite{BM1994}, where AI  limits and a schematic bifurcation curve of $1$-$2$ hole transition of cantori for  the twist map \eqref{2harmonicMap} is presented.

\begin{figure}[htbp]
	\begin{center}
		\includegraphics[width=0.52\textwidth, angle=0]{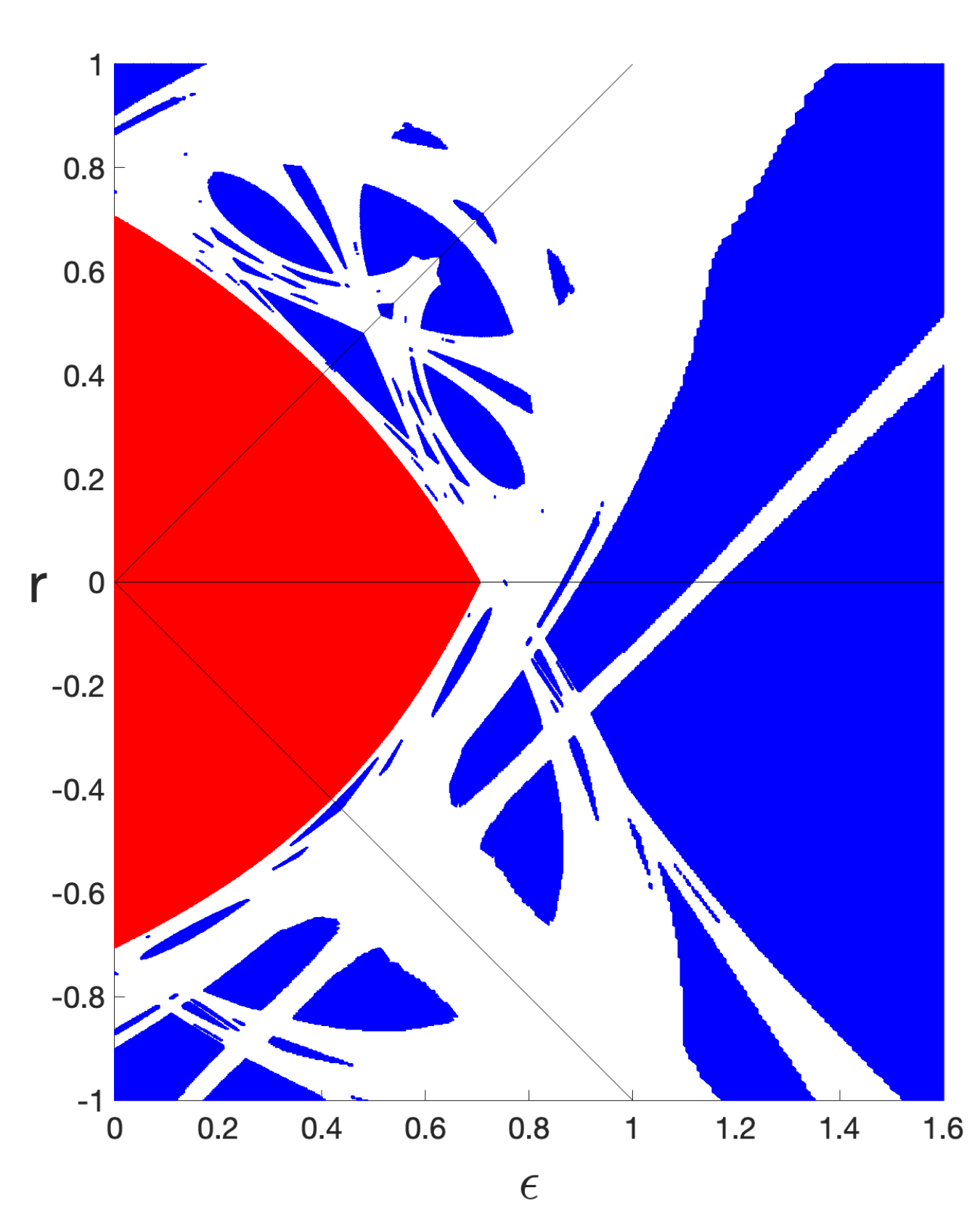}
		\includegraphics[width=0.46\textwidth, angle=0]{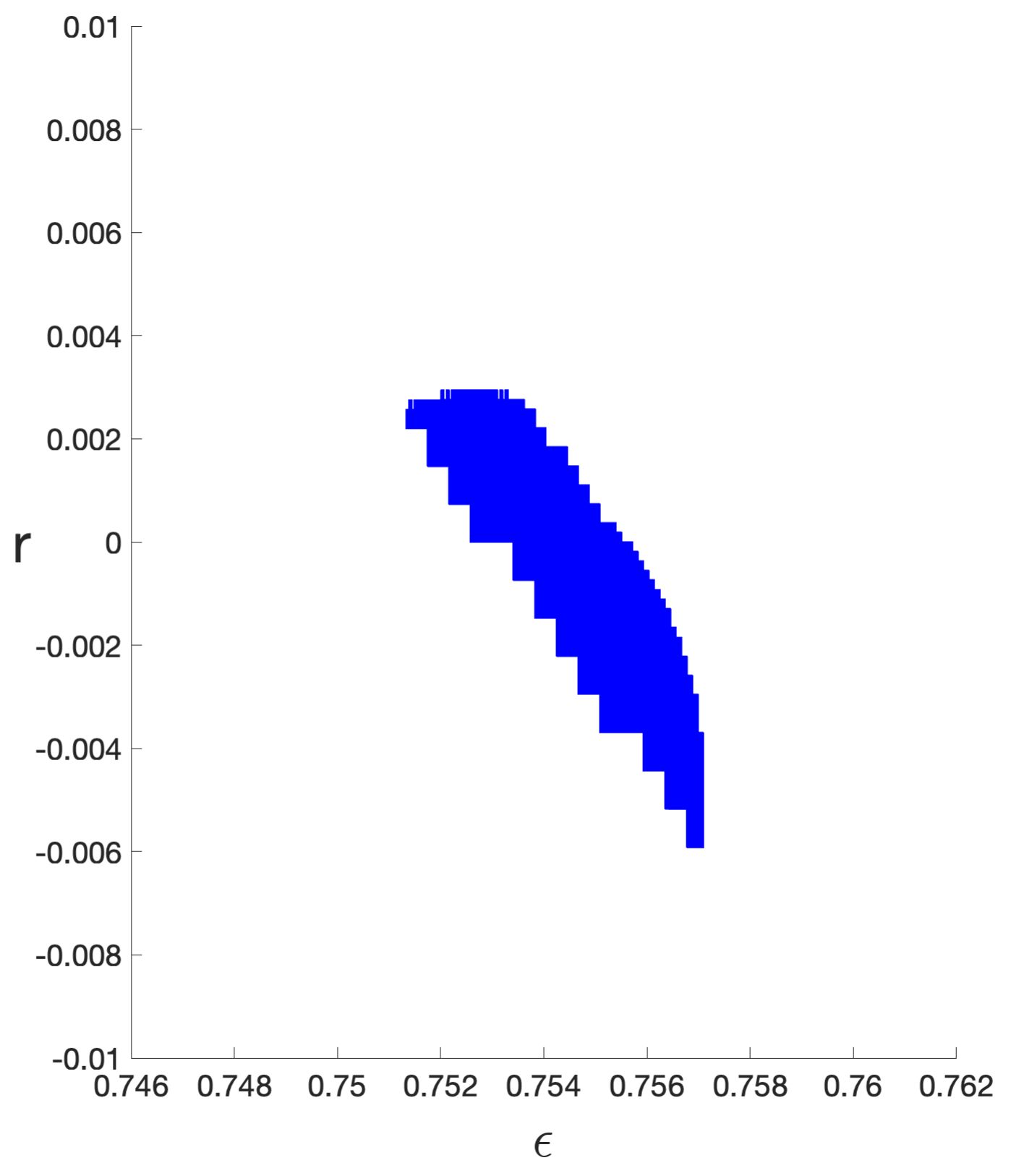}
	\end{center}
	\caption{(Left) Hyperbolic plateaus on  the $(\epsilon, r)$-plane. (Right) Magnification for $0.746\le \epsilon\le 0.762$ and $-0.010\le r\le 0.010$. When $r=0$, the blue part is contained in the attracting period-$3$ window of the quadratic map $Q_a$.}
	\label{epsilon_r_plateaus}
\end{figure}

\paragraph{Plateaus on $(\epsilon, r)$-plane.}

Figure \ref{epsilon_r_plateaus} shows the plateaus on the $(\epsilon, r)$-plane. 
The coloured regions are the hyperbolic plateaus, in which the horseshoe component is in red. The lines $r=\pm \epsilon$  represents $b=\pm 1$ in the parameter $(a,b)$-plane.
 This figure has similar advantage to Figure \ref{semi_disc}{\textemdash}it ``blows up" the  point infinity 
(i.e.  the ``north pole" of the  sphere in Figure \ref{fig:plateaus_sphere} or the point $(1,0)$  in the unit disc in Figure \ref{fig:plateaus_disc})
 into a line (i.e. the $r$-axis). This figure 
 is particularly suitable for features of the parameters near the AI limit.
Notice that the plateaus are symmetric with respect to the lines $r=\epsilon$ and $r=-\epsilon$.

These symmetries are expected from the symmetric roles played by $\epsilon$ and $r$ in the recurrence relation 
\eqref{HenonRecurrence}. These symmetries also indicate that the dynamics for parameters 
near and at the $\epsilon$-axis and $r$-axis are topologically conjugate. More precisely, define $f_{\epsilon, r}(x,y):=H_{1/\epsilon^2, r/\epsilon}(x,y)$, then $f_{\epsilon, r}^{-1}$ is topologically conjugate to $f_{r,\epsilon}$.
 In particular, the H\'{e}non map defined as $f_{\epsilon, r}$ is a hyperbolic horseshoe when $-1/\sqrt{2}<r<1/\sqrt{2}$ and $\epsilon>0$ sufficiently small, also when $0<\epsilon<1/\sqrt{2}$ and $|r|$ sufficiently small. The AI limit, namely $\epsilon=r=0$, is just a special and the simplest case.

Furthermore, these symmetries also indicate that all mathematical analyses for perturbations of $H_{a,b}$ around the limit $b\to 0$ (with parameter $a$ fixed) and around the limit $a\to\infty$ with $b^2/a=constant$ are actually the same:
Taking the limit $a\to \infty$ with $b^2/a=\hat{r}^2$ for $H_{a,b}^{-1}$ means taking the limit $\epsilon\to 0$ for $H_{1/\epsilon^2, \hat{r}/\epsilon}^{-1}$. The latter is equivalent to taking the limit $\epsilon\to 0$ for $H_{1/\hat{r}^2, \epsilon/\hat{r}}$ because $H_{1/\epsilon^2, r/\epsilon}^{-1}$ is topologically conjugate to $H_{1/r^2, \epsilon/r}$. For a fixed nonzero $\hat{r}$, the Jacobian of $H_{1/\hat{r}^2, \epsilon/\hat{r}}$ is $\epsilon/\hat{r}$. Therefore, the limit $\epsilon\to 0$ means the Jacobian tends to zero. In view of (ii) of Theorem \ref{thm:AI2}, we infer that the restriction of $H_{1/\hat{r}^2, \epsilon/\hat{r}}$ to 
\[ \mathcal{S}_\epsilon:=\bigcup_{{\bf x}^\dag\in\Sigma_{0,\hat{r}}}(x_0^*, x_1^*)
\]
 is topologically  conjugate to $\widehat{Q}_{1/\hat{r}^2}$ on $\underleftarrow{\lim}(\Lambda_{1/\hat{r}^2}, Q_{1/\hat{r}^2})$ if  $0<\epsilon<\epsilon_0$ and $1/\hat{r}^2\in\mathsf{Hyp}$. Hence, we conclude

\begin{cor}[Perturbation of quadratic maps to H\'{e}non maps]
 Suppose $a\in\mathsf{Hyp}$, then the restriction of $H_{a,b}$ to $S_{\sqrt{a}b}$ is topologically conjugate to $\widehat{Q}_a$ on $\underleftarrow{\lim}(\Lambda_a, Q_a)$ provided $0<b<\epsilon_0/\sqrt{a}$.
\end{cor}

We end this section with a remark that  for a parameter $a$ with which $Q_a$ has an attracting periodic orbit and sufficiently small $b$, Barge and Holte \cite{BH1995} showed that $H_{a,b}$, when restricted to its attracting set, is topologically conjugate to the shift homeomorphism on the inverse limit space of the interval $[\sqrt{a}(1-a), \sqrt{a}]$ with bonding map $Q_a$.

\section{Proof of theorem \ref{thm:AI2}} \label{sec:proof}

The proof is based on the implicit function theorem.

\subsubsection*{Invertibility of $DF({\bf x}^\dag; 0)$:}

\begin{la} \label{keylemma}
 The linear operator $DF({\bf x}^\dag; 0)$ is invertible with bounded inverse for all ${\bf x}^\dag\in\Sigma_{0, \hat{r}}$.
\end{la}
  \proof
For a given ${\bf x}^\dag\in \Sigma_{0, \hat{r}}$,   the operator $DF({\bf x}^\dag; 0)$ is invertible if and only if 
\begin{equation} 
   \hat{r} \xi_{i-1}+2x_i^\dag\xi_i=\eta_i  \label{non-h} 
\end{equation}
possesses a unique solution $\boldsymbol{\xi}=(\xi_i)_{i\in\Z}\in l_\infty(\Z, \R)$ for any given $\boldsymbol{\eta}=(\eta_i)_{i\in\Z}\in l_\infty(\mathbb{Z},\mathbb{R})$. If  $\hat{r}\not=0$, equation (\ref{non-h}) has a solution 
\begin{eqnarray}
  \xi_i &=& \sum_{N\ge 0} (-\hat{r})^{N}\left(\prod_{k=0}^N \left(2x_{i-k}^\dag\right)^{-1}\right)\eta_{i-N} \label{11oneone} \\
        &=& - \hat{r}^{-1}\sum_{N\ge 0}\left(\left(Q_{1/\hat{r}^2}^{N+1}\right)'(x_{i-N}^\dag)\right)^{-1}\eta_{i-N}, \label{solnon-h}
\end{eqnarray}
which depends on both ${\bf x}^\dag$ and $\boldsymbol{\eta}$. The solution $\xi_i$ is bounded for every $i\in\Z$ because the right hand side of the equality of (\ref{solnon-h}) can be bounded by a geometric series coming from  the expanding property that $|(Q^n_{1/\hat{r}^2})^\prime (x)|\ge C\lambda^n$ for some  constants $C>0$ and 
$\lambda>1$ for all $n\ge 1$ and  all $x\in\Lambda_{1/\hat{r}^2}$.  The  expanding property also guarantees the norm of $DF({\bf x}^\dag; 0)^{-1}$ is uniformly bounded in ${\bf x}^\dag$ since $\|DF({\bf x}^\dag; 0)^{-1}\|$ is equal to $\sup_{i\in\Z}|\xi_i|/\sup_{\boldsymbol{\eta}\in l_\infty(\Z,\R)}\|\boldsymbol{\eta}\|$.  Moreover, the same expanding property implies the solution $(\xi_i)_{i\in\Z}$ is unique: Any homogeneous solution of \eqref{non-h} satisfies 
\[
  \xi_{i-N} = \left( Q_{1/\hat{r}^2}^{N}\right)'(x_{i-N+1}^\dag)~\xi_i
\]
for any positive integer $N$, thus $\xi_{i-N}$ can be arbitrarily large for large $N$ if $\xi_i$ is not zero.

If $\hat{r}=0$,  the operator $D_{\bf x}F({\bf x}^\dag, 0)$ is invertible if and only if 
\[
   2x_i^\dag\xi_i=\eta_i
\]
possesses a unique bounded solution for any given $(\eta_i)_{i\in\Z}\in l_\infty(\mathbb{Z},\mathbb{R})$. It is clear that $\xi_i=\eta_i/ (2x_i^\dag)$ for every integer $i$ since $x_i^\dag$ is either $-1$ or $1$.   
\qed

\subsubsection*{Continuation:}

 If $\Lambda_{1/\hat{r}^2}$ is a compact hyperbolic invariant set for $Q_{1/\hat{r}^2}$, then the linear operator $DF({\bf x}^\dag;  0)$ is invertible for every ${\bf x}^\dag\in \Sigma_{0, \hat{r}}$  by  Lemma \ref{keylemma}.  Since $F({\bf x}^\dag; 0)=0$, by virtue of  the implicit function theorem (e.g. \cite{BP1970}),  there exists $\epsilon_1=\epsilon_1({\bf x}^\dag)$ and a $C^1$ function 
  \[ {\bf x}^*(\cdot;{\bf x}^\dag):\R\to l_\infty(\Z,\R), \qquad \epsilon\mapsto   {\bf x}^*(\epsilon;{\bf x}^\dag)=({\bf x}^*(\epsilon;{\bf x}^\dag)_i)_{i\in\Z}=(x_i^*)_{i\in\Z}
  \]
  such that $F({\bf x}^*(\epsilon; {\bf x}^\dag);\epsilon)=0$ and ${\bf x}^*(0;{\bf x}^\dag)={\bf x}^\dag$ provided $0\le \epsilon<\epsilon_1$. The construction of $F$ tells that ${\bf x}^*(\epsilon;{\bf x}^\dag)$ is an orbit of the map $H_{1/\epsilon^2, r/\epsilon}$.

\subsubsection*{Bijectivity:}

We show first that 
 $\inf_{            {\bf x}^\dag\in\Sigma_{0, \hat{r}}
                 } \epsilon_1({\bf x}^\dag
                         )> 0$.
   There are positive $\delta_1$ and $\epsilon_2$ such that for ${\bf x}$ in the closed ball $\bar{B}({\bf x}^\dag, \delta_1)$ of radius $\delta_1$ centred at ${\bf x}^\dag$ and for $0<\epsilon\le\epsilon_2$, we have
\begin{eqnarray*}
  \|DF({\bf x} ;\epsilon)-DF({\bf x}^\dag;0)\|
&=& \sup_{\|\boldsymbol{\xi}\|=1} \left|\epsilon\xi_{i+1}+2(x_i -x_i^\dag)\xi_i+(r(\epsilon)-\hat{r})\xi_{i-1}\right| \\
&\le & 1/ (2~ \|DF({\bf x}^\dag; 0)^{-1}\|) 
\end{eqnarray*}
and
\begin{equation}
\|F({\bf x}^\dag; \epsilon)\| <\delta_1 /(2~ \|DF({\bf x}^\dag; 0)^{-1}\|). \label{etadelta}
\end{equation}
For each ${\bf x}^\dag\in\Sigma_{0,\hat{r}}$ and $\epsilon$, define a map $\mathcal{F}(\cdot;{\bf x}^\dag,\epsilon): l_\infty(\Z,\R) \to l_\infty(\Z, \R)$, ${\bf x}\mapsto {\bf x}-DF({\bf x}^\dag; 0)^{-1}F({\bf x}; \epsilon)$. Thus, for ${\bf v}$, ${\bf w}\in \bar{B}({\bf x}^\dag,\delta_1)$ and $0<\epsilon\le\epsilon_2$, we get 
\begin{eqnarray}
 \lefteqn{ \left\|(\mathcal{F}({\bf v};{\bf x}^\dag,\epsilon)-\mathcal{F}({\bf w};{\bf x}^\dag,\epsilon)\right\|  } \nonumber \\
&\le &  \sup_{\bar{\bf u}\in \{{\bf w}+t({\bf v}-{\bf w}):~ 0\le t\le 1\}}\left\|DF({\bf x}^\dag; 0)^{-1} \left(
                          DF({\bf x}^\dag; 0)-  DF(\bar{\bf u}; \epsilon)                                                                                   \right)\right\|\cdot\|({\bf v}-{\bf w})\|\nonumber\\
&\le& \|{\bf v}-{\bf w}\| /2 \label{401}
\end{eqnarray}
and then
\begin{eqnarray}
 \left\|\mathcal{F}({\bf x};{\bf x}^\dag,\epsilon)-{\bf x}^\dag\right\|
&\le & \left\| \mathcal{F}({\bf x};{\bf x}^\dag,\epsilon)-\mathcal{F}({\bf x}^\dag;{\bf x}^\dag,\epsilon)\right\|+\left\|\mathcal{F}({\bf x}^\dag;{\bf x}^\dag,\epsilon)-{\bf x}^\dag\right\|\nonumber\\
&\le& \|{\bf x}-{\bf x}^\dag\| /2+\|DF({\bf x}^\dag; 0)^{-1}\| ~\|F({\bf x}^\dag; \epsilon)\|\label{402}\\
&<&\delta_1 \nonumber
\end{eqnarray}
for any ${\bf x} \in \bar{B}({\bf x}^\dag,\delta_1)$.
This implies that $\mathcal{F}(\cdot;{\bf x}^\dag,\epsilon)$ is a contraction map with contraction constant at least $1/2$ on $\bar{B}({\bf x}^\dag,\delta_1)$ for any ${\bf x}^\dag\in\Sigma_{0,\hat{r}}$
 and $0\le \epsilon\le\epsilon_2$. Hence, $\inf_{{\bf x}^\dag\in\Sigma_{0, \hat{r}}}\epsilon_1({\bf x}^\dag)\ge\epsilon_2>0$.

Secondly, the radius $\delta_1$ is independent of ${\bf x}^\dag$ and of $\epsilon$, and ${\bf x}^*(\epsilon; {\bf x}^\dag)$  is the unique fixed point in $\bar{B}({\bf x}^\dag,\delta_1)$ for $\mathcal{F}(\cdot;{\bf x}^\dag,\epsilon)$.  

Thirdly, notice that $\Lambda_{1/\hat{r}^2}$ is expansive. That is, there exists $\tau>0$ such that for  distinct $x$, $\tilde{x}\in\Lambda_{1/\hat{r}^2}$, there is some $n\ge 0$ such that $|Q_{1/\hat{r}^2}^n(x)-Q_{1/\hat{r}^2}^n(\tilde{x})|>\tau$. This implies that $\Sigma_{0, \hat{r}}$ is uniformly discrete.
  (A subset $\Sigma$ of $l_\infty(\Z, \R)$ is {\it uniformly discrete} means that there exists $\tau^\prime>0$ such that whenever ${\bf u}$ and ${\bf v}$ are distinct points of $\Sigma$, then $\|{\bf u}-{\bf v}\| >\tau^\prime$.)  Because $\Sigma_{0,\hat{r}}$ is uniformly discrete, the balls $\bar{B}({\bf x}^\dag,\delta_1)$, ${\bf x}^\dag\in\Sigma_{0,\hat{r}}$, are disjoint in $l_\infty (\Z, \R)$ provided that $\delta_1$ is sufficiently small (by taking smaller $\epsilon_2$ if necessary). It follows that the mapping 
  \[ \Psi_\epsilon :     \Sigma_{0, \hat{r}} \to  \bigcup_{
                                                                {\bf x}^\dag\in\Sigma_{0, \hat{r}}
                                                            } {\bf x}^*(\epsilon; {\bf x}^\dag)
   \]
is a bijection provided $0\le \epsilon<\epsilon_2$.

The above proves the statement (i) of  the theorem for $|x_i^*-x_i^\dag|<\delta_1$ and $0\le \epsilon<\epsilon_2$.

\subsubsection*{Topological conjugacy:}

  Let the projection ${\bf x}=(\cdots,x_{-1},x_0,x_1,\cdots)\mapsto (x_0, x_{-1}) \in\mathbb{R}^{2}$ be  denoted by $\pi_0$.  
 We assert in the next paragraph that   the  composition of mappings  
  \[ {\bf x}^\dag\ \stackrel{\Psi_\epsilon}{\longmapsto}\ {\bf x}^*(\epsilon;{\bf x}^\dag)\ \stackrel{\pi_0}{\longmapsto}~ (x^*_0, x^*_{-1}) 
  \]
  is a homeomorphism with the product topology, and that  the  following diagram 
\[
  \begin{matrix}
        \Sigma_{0, \hat{r}}& \mapright{\sigma} &   \Sigma_{0, \hat{r}} \cr
         \mapldown{\pi_0\circ\Psi_\epsilon}& &\maprdown{\pi_0\circ\Psi_\epsilon} \cr
             \mathcal{A}_{\epsilon}& \mapright{H_{1/\epsilon^2,r/\epsilon}}& \mathcal{A}_{\epsilon}
  \end{matrix}
\]        
commutes when $0<\epsilon<\epsilon_0$ for some $\epsilon_0\le \epsilon_2$ and $\delta_0=\delta_1/2$, where the set $\mathcal{A}_{\epsilon}$ is compact, invariant under $H_{1/\epsilon^2,r/\epsilon}$, and is defined by 
  \[ \mathcal{A}_{\epsilon}:=\bigcup_{
                                                {\bf x}^\dag\in   \Sigma_{0, \hat{r}}
                                            }                                       \pi_0\circ {\bf x}^*(\epsilon;{\bf x}^\dag).
\]
In other words, $\pi_0\circ\Psi_\epsilon$ acts as the topological conjugacy. (Hence, the statements (i) and (ii) of Theorem \ref{thm:AI2} are proved.)

  The projection  $\pi_0$ certainly is continuous. It is bijective since the sequence ${\bf x}^*(\epsilon;{\bf x}^\dag)$ is an orbit of the diffeomorphism $H_{1/\epsilon^2, r/\epsilon}$ and is uniquely determined by the initial point $(x_0^*, x_{-1}^*)$. We showed that  $\Psi_\epsilon$ is bijective. Note that the compactness of $\Lambda_{1/\hat{r}^2}$ implies the compactness of  $\Sigma_{0, \hat{r}}$ in the product topology. We shall show in Lemma \ref{lala} below that  $\Psi_\epsilon $ is continuous when $0<\epsilon<\epsilon_0$. As a consequence,  $\pi_0\circ\Psi_\epsilon$ is a continuous bijection from a compact space to a Hausdorff space, thus a homeomorphism. Note also that $Q_{1/\hat{r}^2}$-invariance of $\Lambda_{1/\hat{r}^2}$ implies 
  $\sigma$-invariance of  $\Sigma_{0, \hat{r}}$.
Now, 
\begin{eqnarray*}
F_{i+1}({\bf x}; \epsilon)
 &=&  \epsilon x_{i+2}-(1-x_{i+1}^2)+r x_{i}\\
            &=&\epsilon\sigma({\bf x})_{i+1}-(1-\sigma ({\bf x})_i^2)+r\sigma({\bf x})_{i-1} \\
          &=& F_i(\sigma ({\bf x}); \epsilon).
\end{eqnarray*}
This means that 
\begin{equation}
  \sigma\circ F({\bf x}; \epsilon)=F(\sigma({\bf x}); \epsilon). \label{=23}
  \end{equation} 
  We showed that $F(\cdot; \epsilon)$ has a unique zero at ${\bf x}^*(\epsilon; {\bf x}^\dag)$ in $\bar{B}({\bf x}^\dag, \delta_1)$, so does $F(\sigma(\cdot);\epsilon)$ by  \eqref{=23}. This implies that $F(\cdot; \epsilon)$ has a unique zero at $\sigma ( {\bf x}^*(\epsilon; {\bf x}^\dag))$ in $\bar{B}(\sigma({\bf x}^\dag),\delta_1)$. However, $F(\cdot; \epsilon)$ has been shown to have a unique zero at ${\bf x}^*(\epsilon; \sigma({\bf x}^\dag))$ in $\bar{B}(\sigma({\bf x}^\dag), \delta_1)$. Thence, $\sigma ({\bf x}^*(\epsilon; {\bf x}^\dag))$ must equal ${\bf x}^*(\epsilon; \sigma({\bf x}^\dag))$, and  the commutativity of the diagram follows immediately.

\begin{la} \label{lala}
 There exists $\epsilon_0\le \epsilon_2$ such that the  map  $\Psi_\epsilon$ is continuous in the product topology  for $0<\epsilon<\epsilon_0$. 
\end{la}
\proof
For $0<\epsilon\le \epsilon_2$, the  closed ball $\bar{B}({\bf x}^\dag, \delta_1)=(\bar{B}({x}_i^\dag, \delta_1))_{i\in\Z}$ contains ${\bf x}^*(\epsilon; {\bf x}^\dag)=({\bf x}^*(\epsilon; {\bf x}^\dag)_i)_{i\in\Z}$ for each ${\bf x}^\dag$ in $l_\infty(\Z, \R)$, and if $F({\bf x};\epsilon)=0$ then ${\bf x}={\bf x}^* (\epsilon; {\bf x}^\dag)$ whenever ${\bf x}\in \bar{B}({\bf x}^\dag, \delta_1)$. 
Let 
\[ \bigcap_{i=-n}^n  H_{1/\epsilon^2, r/ \epsilon}^{-i} \left(
                            \bar{B}({x}_i^\dag, \delta_1) \right)=:\bar{B}( {x}_0^\dag, \delta_1)^{n}.
\]
Then, we have a nested sequence of closed sets
\[
 \bar{B}( {x}_0^\dag, \delta_1)\supset \bar{B}( {x}_0^\dag, \delta_1)^{1}\supset \bar{B}( {x}_0^\dag, \delta_1)^{2} \supset\cdots. 
  \]
  In particular,  $\lim_{n\to\infty}\bar{B}( { x}_0^\dag, \delta_1)^{n}=:\bar{B}( { x}_0^\dag, \delta_1)^{\infty}$ consists of only one point, which is ${\bf x}^*(\epsilon; {\bf x}^\dag)_0$. This is because  an orbit ${\bf x}$ of $H_{1/\epsilon^2, r/\epsilon}$ with initial point lying in $\bar{B}( {x}_0^\dag, \delta_1)^{\infty}$ must have $x_i\in \bar{B}( {x}_i^\dag, \delta_1)$ for every integer $i$, but by the uniqueness, it must be ${\bf x}^*(\epsilon;{\bf x}^\dag)_0$. 
  
There exists $\epsilon_0\le\epsilon_2$ such that the statement in the above paragraph  remains true if the radius of ball centred at ${\bf x}^\dag$ is replaced by $\delta_0=\delta_1/2$ and $0<\epsilon<\epsilon_0$.
  
Now, assume $( {\bf x}^{(m)} )_{m>0}$ is a convergent sequence to ${\bf x}^\dag$ in $\Sigma_{0,\hat{r}}$ with the product topology, then 
 there is an integer $N=N(m)>1$ for $m$ large enough such that $|x_i^{(m)}-x_i^\dag|<\delta_0$  for $|i|\le N$. Certainly, 
  $  \bigcap_{i=-\infty}^\infty  H_{1/\epsilon^2, r/ \epsilon}^{-i} \left(
                            \bar{B}({ x}_i^{(m)}, \delta_0) \right)={\bf x}^*(\epsilon;{\bf x}^{(m)})_0$ if $0<\epsilon<\epsilon_0$. 
  Notice that $\bar{B}(x_i^{(m)}, \delta_0)\subset \bar{B}(x_i^\dag, \delta_1)$ for $|i|\le N(m)$. Thus, for $0<\epsilon<\epsilon_0$, we get $\bar{B}(x_0^{(m)}, \delta_0)^{N(m)}\subset \bar{B}(x_0^\dag, \delta_1)^{N(m)}$.
Because ${\bf x}^*(\epsilon;{\bf x}^{(m)})_0$ is contained in $\bar{B}(x_0^{(m)}, \delta_0)^{N(m)}$, it is arbitrarily close to ${\bf x}^*(\epsilon;{\bf x}^\dag)_0$ for sufficiently large $m$. This shows the mapping ${\bf x}^\dag\mapsto {\bf x}^*(\epsilon;{\bf x}^\dag)_i$ is continuous when $i=0$. Similar procedure shows the mapping is continuous for each integer $i$. Hence $\Psi_\epsilon$ is continuous for $0<\epsilon<\epsilon_0$. 
 \qed\\

\subsubsection*{Hausdorff convergence:}

The statement (iii) is straightforward since ${\bf x}^*(\epsilon; {\bf x}^\dag)\to {\bf x}^\dag$ in $l_\infty(\Z, \R)$ when $\epsilon\to 0$.

The proof of Theorem \ref{thm:AI2} is completed.
\qed

\section*{\small ACKNOWLEDGMENTS}

Z. Arai was partially supported by JSPS KAKENHI grant numbers JP25K00919, JP23K17657, JP19KK0068. Y.-C. Chen was partially supported by MOST grant number 112-2115-M-001-005. Y.-C. Chen  thanks Yongluo Cao  for useful conversations  about perturbing the quadratic map to the H\'enon map, and thanks for the hospitality of Hokkaido and Soochow Universities during his visits.

\section*{\small AUTHOR DECLARATIONS}

The authors have no conflicts to disclose.

\section*{\small DATA AVAILABILITY}

Data sharing is not applicable to this article as no new data were created or analyzed in this study.


\begin{thebibliography}{99}

\bibitem{Arai2007}
Z. Arai,
On hyperbolic plateaus of the H\'{e}non map,  
{\it Experiment. Math.}  {\bf 16} (2007) 181--188. 

\bibitem{Aubr1995}
S. Aubry, 
Anti-integrability in dynamical and variational problems,
{\it Physica D} {\bf 86} (1995) 284--296.

\bibitem{AA1990}
S. Aubry, G. Abramovici, 
Chaotic trajectories in the standard map: the concept of  anti-integrability,
{\it Physica D} {\bf 43} (1990) 199--219.

\bibitem{BCM2013}
C. Baesens, Y.-C. Chen, R.S. MacKay, 
Abrupt bifurcations in chaotic scattering: view from the anti-integrable limit,
{\it Nonlinearity} {\bf 26} (2013) 2703--2730.

\bibitem{BM1993}
C. Baesens, R.S. MacKay,
Cantori for multiharmonic maps,
{\it Phys. D} {\bf 69} (1993) 59--76. 

\bibitem{BM1994}
C. Baesens, R.S. MacKay, 
The one to two-hole transition for cantori, 
{\it Phys. D} {\bf 71} (1994) 372--389. 

\bibitem{BH1995}
M. Barge, S. Holte, 
Nearly one-dimensional H\'{e}non attractors and inverse limits,
{\it Nonlinearity} {\bf 8} (1995) 29--42.

\bibitem{BM1997}
S.V. Bolotin, R.S. MacKay, 
Multibump orbits near the anti-integrable limit for Lagrangian systems,
{\it Nonlinearity} {\bf 10} (1997) 1015--1029.

\bibitem{BM2000}
S.V. Bolotin, R.S. MacKay,
Periodic and chaotic trajectories of the second species for the $n$-centre problem,
{\it Celestial Mechanics and Dynamical Astronomy} {\bf 77} (2000)  49--75.

\bibitem{BT2015}
S.V. Bolotin, D.V. Treschev, 
The anti-integrable limit,
{\it Russian Mathematical Surveys} {\bf 70} (2015) 975--1030.

\bibitem{Bowe1970}
R. Bowen,
Topological entropy and axiom A, {\it Proceedings of Symposia in Pure Mathematics},  American Mathematical Society, Providence, RI, 1970, pp. 23-41.

\bibitem{BP1970}
A.L. Brown, A. Page, 
{\it Elements of Functional Analysis}, 
Van Nostrand Reinhold Company, London, 1970.

\bibitem{Chen2003}
Y.-C. Chen, 
Multibump orbits continued from the anti-integrable limit for Lagrangian systems.
{\it  Regul. Chaotic Dyn.} {\bf 8} (2003) 243--257.

\bibitem{Chen2004}
Y.-C. Chen, 
Anti-integrability in scattering billiards, 
{\it Dynamical Systems} {\bf 19} (2004) 145--159.

\bibitem{Chen2005}
Y.-C. Chen, 
Bernoulli shift for second order recurrence relations near the anti-integrable limit,
{\it Discrete and Continuous Dynamical Systems B} {\bf 5} (2005)  587--598.

\bibitem{Chen2006}
Y.-C. Chen,
Smale horseshoe via the anti-integrability,
 {\it Chaos, Solitons \& Fractals} {\bf 28} (2006) 377--385.

\bibitem{Chen2007}
Y.-C. Chen, 
Anti-integrability for the logistic maps, 
{\it Chinese Ann. Math. B} {\bf 28} (2007) 219--224.

\bibitem{Chen2010}
Y.-C. Chen,
On topological entropy of billiard tables with small inner scatterers,
{\it Advances in Mathematics}  {\bf 224} (2010) 432--460

\bibitem{Chen2018}
Y.-C. Chen, 
A proof of Devaney-Nitecki region for the H\'{e}non mapping using the anti-integrable limit, 
{\it Advances in Dynamical Systems and Applications} {\bf 13} (2018) 33--43.

\bibitem{CFS1977}
R.C. Churchill, J. Franke, J. Selgrade,
A geometric criterion for hyperbolicity of flows,
{\it Proc. Amer. Math. Soc.} {\bf 62} (1977) 137--143.

\bibitem{dMvS1993}
W. de Melo,  S. van Strien, 
{\it One-dimensional Dynamics}, Springer-Verlag, 1993.

\bibitem{DN1979}
R.L. Devaney, Z. Nitecki, 
Shift automorphisms in the H\'{e}non mapping, 
{\it Comm. Math. Phys.} {\bf 67} (1979) 137--146.

\bibitem{GS1997}
J. Graczyk, G. \'Swiatek, 
Generic hyperbolicity in the logistic family, 
 {\it Ann. of Math.} (2) {\bf 146} (1997)  1--52. 

\bibitem{Ishi2008}
Y. Ishii,
Hyperbolic polynomial diffeomorphisms of $\mathbb{C}^2$. I. A non-planar map,
{\it Advances in Mathematics} {\bf 218} (2008)  417--464.

\bibitem{KH1995}
A. Katok, B. Hasselblatt, 
{\it Introduction to the Modern Theory of Dynamical Systems}, Cambridge University Press, Cambridge, 1995.

\bibitem{MM1992}
R.S. MacKay, J.D. Meiss,
Cantori for symplectic maps near the anti-integrable limit, 
{\it Nonlinearity} {\bf 5} (1992) 149--160.

\bibitem{Mane1985}
R. Ma\~{n}\'{e}, 
Hyperbolicity, sinks and measure in one-dimensional dynamics, 
{\it Comm. Math. Phys.} {\bf 100} (1985) 495--524. 

\bibitem{Misi1981}
M. Misiurewicz, 
Absolutely continuous measures for certain maps of an interval, 
{\it Inst. Hautes Etudes Sci. Publ. Math.} {\bf 53} (1981) 17--51. 

\bibitem{MNTU2000}
S. Morosawa, Y. Nishimura, M. Taniguchi, T. Ueda, 
{\it Holomorphic Dynamics},  Cambridge University Press, Cambridge, 2000. 

\bibitem{Mumm2008}
P.P. Mummert, 
Holomorphic shadowing for H\'{e}non maps,
{\it Nonlinearity} {\bf 21} (2008) 2887--2898. 

\bibitem{Qin2001}
W.-X. Qin,
Chaotic invariant sets of high-dimensional H\'{e}non-like maps,
{\it J. Math. Anal. Appl.} {\bf 264} (2001)  76--84.

\bibitem{SS1974}
R.J. Sacker, G.R. Sell,
Existence of dichotomies and invariant splitting s for linear differential systems I,
{\it J. Diff. Eqns.}  {\bf 15} (1974)  429--458.

\bibitem{Sing1978}
D. Singer,
Stable orbits and bifurcation of maps of the interval,  
{\it SIAM J. Appl. Math.} {\bf 35} (1978)  260--267. 

\bibitem{SM1998}
D. Sterling, J.D. Meiss, 
Computing periodic orbits using the anti-integrable limit, 
{\it Phys. Lett. A} {\bf 241} (1998) 46--52.

\end{thebibliography}
\end{document}